\documentclass[11pt]{amsart}
\usepackage{amssymb,latexsym,eufrak,amsmath,amscd,
graphicx,maple2e}
\usepackage{amsthm}

\theoremstyle{plain}
\newtheorem*{Main}{Main Theorem}

\theoremstyle{plain}

\newtheorem{Cor}{Corollary}
\newtheorem{Lem}{Lemma}

\theoremstyle{definition}
\newtheorem{Def}{Definition}
\newtheorem{notation}{Notation}
\newtheorem{rmk}{Remark}

\theoremstyle{definition}
\newtheorem{Ex}{Example}

\newcommand{\floor}[1]{\lfloor #1 \rfloor}
\newcommand{\tensor}{\otimes}
\newcommand{\one}{{\mathbf 1}} 

\begin{document}
	\title{Multiplier Ideals of Monomial Ideals}
	\author{J. A. Howald}
	\date{February 2, 2000}
	\thanks{I would like to thank Robert Lazarsfeld for suggesting this problem,
		and for many valuable discussions.}


	\begin{abstract}
		In this note we discuss a simple algebraic 
		calculation of the multiplier 
		ideal associated to a monomial ideal in affine n-space.
		We indicate how this result allows one to compute
		not only the multiplier ideal but also the 
		log canonical threshold
		of an ideal in terms of its Newton polygon.
	\end{abstract}

	\maketitle
	
	\section*{Introduction} \label{S:intro}

	Multiplier ideals have become quite important in higher
	dimensional geometry, because of their strong vanishing
	properties, (cf 
	\cite{Angehrn-Siu}, \cite{Demailly}, \cite{Demailly-Kollar},
	\cite{Ein}, \cite{Ein-Lazarsfeld}, \cite{Siu}).
	They reflect the singularity
	of a divisor, ideal sheaf, or metric.
	It is however fairly difficult to calculate
	multiplier ideals explicitly, even in the simplest cases: the
	algebraic definition of the multiplier ideal associated to an
	arbitrary ideal sheaf \( {\mathfrak a} \) requires that we
	construct a log resolution of \( {\mathfrak a} \) and perform
	calculations on the resolved space.  In this note, we
	compute the multiplier ideal associated to an arbitrary
	monomial ideal \( {\mathfrak a} \).  Like \( {\mathfrak a} \), it can
	be described in combinatorial and linear-algebraic terms.

	We begin with some definitions.  Let \( X \) be a smooth
	quasiprojective complex algebraic variety.  Let 
	\( {\mathfrak a} \subset {\mathcal O}_{X} \) be any ideal sheaf.  By a {\bf log 
	resolution} of \( {\mathfrak a} \), we mean a proper birational map
	\( f: {Y} \rightarrow {X} \) with the property that \( Y \) is smooth and
	\(f^{-1}({\mathfrak a}) = {\mathcal O}_{Y}(-E) \), where \( E \) is an effective Cartier divisor, 
	and \( E + exc(f)\) has normal crossing support.
	
	\newcommand{\multr}[1]{{\mathcal J}(r \cdot #1)}
	\newcommand{\mult}[1]{{\mathcal J}(#1)}
	\begin{Def} \label {D:MultiplierI}
	    Let \( {\mathfrak a} \subset {\mathcal O}_{X} \) be an ideal sheaf
	    in \( X \), and let \( f: {Y} \rightarrow {X} \) be a log
	    resolution of \( {\mathfrak a} \), with 
	    \(f^{-1}({\mathfrak a}) = {\mathcal O}_{Y}(-E) \).  Let \( r>0 \)
	    be a rational number.  We define the multiplier ideal of
	    \( {\mathfrak a} \) with coefficient \( r \) to be:
	    
	    \[ \multr{{\mathfrak a}}  = f_{*}{\mathcal O}_{Y}(K_{{Y} / {X}}-\floor{rE}).\]		
	    
	    Here \( K_{{Y} / {X}} = K_{Y} -f^{*}K_{X} \) is the relative canonical 
	    bundle, and \( \floor{-} \) is the round-down for \( {\mathbb 
	    Q}-\)divisors.
	    That \( \multr{{\mathfrak a}} \) is an ideal sheaf follows from the 
	    observation that 
	    \( {\mathcal O}_{Y}(K_{{Y} / {X}}-\floor{rE}) \)
	    is a subsheaf of \( {\mathcal O}_{Y}(K_{{Y} / {X}}) \):  since 
	    \( f_{*}({\mathcal O}_{Y}(K_{{Y} / {X}})) = {\mathcal O}_{X} \), 
	    \( \multr{{\mathfrak a}} \subset {\mathcal O}_{X} \).
	    We write 
	    \( \mult{\mathfrak a} \) for 
	    \( \mult{1 \cdot {\mathfrak a}}\).
	\end{Def}
	We will now specialize to the case 
	\( X={\mathbb A}^{n} \).
	
	\begin{Def}
		Let \( {\mathfrak a} \subset {\mathbb C}[x_{1},\cdots,x_{n}] \) be a monomial 
		ideal.  We will regard \( {\mathfrak a} \) as a subset of the lattice 
		\( L = {\mathbb N}^{n}\) of monomials.  The {\bf Newton Polygon} \( P \)
		of \( {\mathfrak a} \) is the convex hull of this subset of \( L \), considered
		as a subset of \( L \tensor {\mathbb R} = {\mathbb R}^{n}\).  It is an 
		unbounded region.  \( P \cap L \) is the set of monomials in the 
		integral closure of the ideal \( {\mathfrak a} \) \cite{Eisenbud}.
	\end{Def}
	
	\begin{notation}
		We write \( \one \) for the vector \( (1,1,\ldots,1) \), which is 
		identified with the monomial \( x_{1}x_{2}\ldots x_{n} \).  The 
		associated divisor \( div(\one) \) is the union of the coordinate 
		axes.  We 
		use Greek letters (\( \lambda \in L \)) for elements of \( L \)  or  
		\( L \tensor {\mathbb R} \), and 
		exponent notation \( x^{\lambda} \) for the associated monomials.
		For any subset \( P \) of \( L \tensor {\mathbb R} \), we define 
		\( rP \) ``pointwise:''
		\[ rP= \{ r\lambda : \lambda \in P \}. \]
		We write \( Int(P) \) for the topological interior of \( P \), and 
		\( \floor{P} \) for \( \{ x^{\floor{\lambda}} : \lambda \in P \} \).
	\end{notation}
	
	We regard the Newton polygon ``officially'' as a subset of the real 
	vector space \( L \tensor {\mathbb R} = {\mathbb R}^{n}\); 
	the interior operation \( Int(P) \) relies on the real topology of 
	this vector space.  However, we don't always carefully 
	distinguish \( P \) from the collection of its lattice points \( P \cap L \),
	or from the collection of their associated monomials \( \{ x^{\lambda} : 
	\lambda \in P \cap L \} \).  
	
	Here is our main result:

	\begin{Main} \label {T:Main}
		Let \( {\mathfrak a} \subset {{\mathcal O}_{{\mathbb A}^{n}}}\) be a monomial 
		ideal.  Let \( P \) be its Newton polygon.  Then \(\multr{{\mathfrak 
		a}} \) is a monomial ideal, and contains exactly the following 
		monomials:
		\[ \multr{{\mathfrak a}}  = \{x^{\lambda} : \lambda + \one \in Int(rP) 
		\cap L \}.\]
	\end{Main}

	\begin{rmk}
		The right hand side, \( \{x^{\lambda} : \lambda + \one \in Int(rP) 
		\cap L \} \), could instead be called \( \floor{rP} \).  We state the
		theorem as we do in order to emphasize the monomial \( \one \), 
		which is independently important.
	\end{rmk}
		
	\begin{Ex} \label {Ex:Principle}
		If \( {\mathfrak a} \) is generated by a single monomial, \( x^{\lambda} \), 
		then the polygon \( P \) 
		is the positive orthant translated upward to \( \lambda \), and  
		\[ \mult{{\mathfrak a}}  = \floor{P} = P = {\mathfrak a}.\]
	\end{Ex}
	This is not surprising, because in this case \( {\mathfrak a} \) is already a divisor 
	with normal crossing support.

	\begin{Ex}
		Let us calculate the multiplier ideal of \( (x^{8},y^{6}) \).  The 
		Newton polygon is pictured in Figure 1.  The distinguished integer vectors 
		\( \lambda \) are those with the property that \( \lambda + \one \in 
		Int(P) \).  From Figure 1, we conclude 
		\[ \mult{x^{8},y^{6}} = (x^{6}, x^{5}y, x^{4}y^{2}, 
		x^{2}y^{3},xy^{4},y^{5}). \]
		Notice that \( x^{3}y^{2} \) is almost but not quite in \( 
		\mult{x^{8},y^{6}} \), because \( x^{4}y^{3} \) lies on the 
		boundary, not the interior, of the Newton polygon.
	\end{Ex}
	
	\begin{figure}
		\hskip 1.4 in

		\begin{maplegroup}
		\begin{center}
		\mapleplot{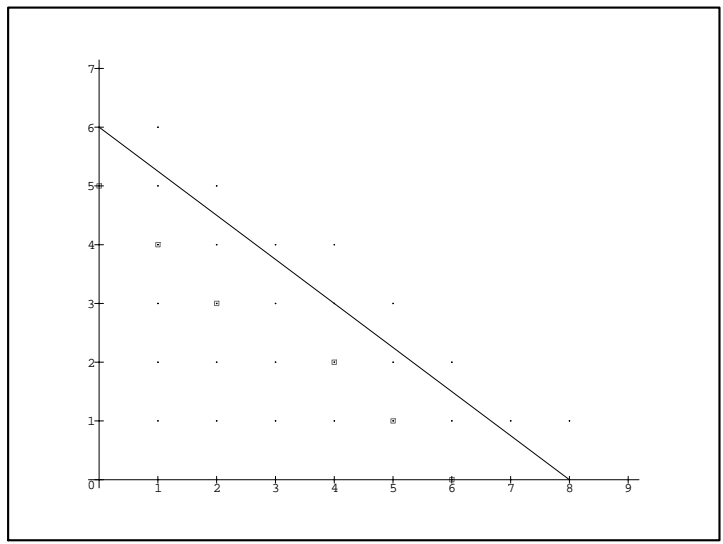}
		\end{center}
		\end{maplegroup}

		\caption{The Multiplier ideal of \( (x^8, y^6) \)}
	\end{figure}
	\begin{Ex} \label {E:EgFrac1}
		Let \( (a_{i})_{i \in n} \) be positive integers, and let \( 
		{\mathfrak a} = (x_{1}^{a_{1}},\ldots,x_{n}^{a_{n}}) \).  One might 
		call this a ``diagonal ideal.''  The only interesting face of 
		the Newton polygon \( P \) of \( {\mathfrak a} \) is defined by a single dual vector 
		\( v=(\frac{1}{a_{1}},\ldots,\frac{1}{a_{n}}) \).  
		Therefore \( \mult{{\mathfrak a}} \) contains the monomials \( 
		\{x^{\lambda} : v \cdot (\lambda + \one) >1 \} \).
		See \cite[example 5.10]{Demailly}, for an analytic perspective on this same result.
		In this expression, the term 
		\( v \cdot \one \) ( \(=\frac{1}{a_{1}}+\cdots+\frac{1}{a_{n}}\)) 
		may be familiar: It is the log-canonical threshold of \( {\mathfrak a} \)
		(see below).
		
	\end{Ex}

	\begin{Ex}
		Let \( g \in {\mathcal O}_{{\mathbb A}^{n}} \) be an arbitrary 
		polynomial.  One might hope that the multiplier ideal associated
		to the (non-monomial) ideal \( (g) \) would be 
		identical to that associated to
		the monomial ideal \({\mathfrak a}_{g}\) generated by the monomials
		appearing in \( g \).  This is not true.  Consider 
		\( g=(x+y)^{n} \) in \( {\mathbb C}[x,y] \).  By a linear change of 
		coordinates in which \( z=x+y \) we obtain \( g=z^{n} \), and can
		calculate \( \mult{(g)} \) in terms of \( z \).  This gives 
		\( \mult{(g)} = (g) \neq \mult{{\mathfrak a}_{g}} \).  
		
		Notice however that for any polynomial \( g \), \( (g) \subset {\mathfrak a}_{g}\).  
		It is not 
		difficult to show that 
		\( \multr{(g)} \subset \multr{{\mathfrak a}_{g}} \) for all \( r \).
		This containment is almost always strict, but it does become an 
		equality if both
		\( r < 1 \) and the coefficients of \( g \) are sufficiently general.
		
		These conditions guarantee that the multiplicity of the		
		\( {\mathbb Q}-\)divisor \( r \cdot (g=0) \) is less than
		one away from the zeroes of 
		\( {\mathfrak a}_{g} \).
		
	\end{Ex}
	
	\begin{Ex} \label {E:lct}
	    Let \( {\mathfrak a} \) be a monomial ideal in \( {\mathbb A}^{n}
	    \), and let \( P \) be its Newton polygon.  The log
	    canonical threshold \( t \) of \( {\mathfrak a} \) is defined
	    to be 
	    \[ t = sup \{r: \multr{{\mathfrak a}} \neq {\mathcal O}_{X} \}.\]
	    See \cite{Kollar} or \cite{Demailly-Kollar} for a detailed discussion of this concept.
	    The Main Theorem shows that this must be equal to \( sup
	    \{r: \one \notin rP \} \) (provided that \( \multr{{\mathfrak a}}\)
	    is nontrivial--the trivial case is an annoying exception).  
	    Thus the log canonical
	    threshold is the reciprocal of the (unique) number \( m \)
	    such that the boundary of \( P \) intersects the diagonal
	    in the \( {\mathbb R}^{n} \) at the point \( m \one \).  In other words,
	    in order to calculate the threshold, we need only find
	    where \( P \) intersects the diagonal.  Arnold calls this
	    number \( m \) of the intersection point the
	    ``remoteness'' of the polygon.  In \cite{Arnold}, he proves
		that \( m = \frac{1}{t} \),
		in order to analyze asymptotic oscillatory integrals.
	    
	\end{Ex}

	\begin{Ex} \label {E:EgFrac2}
	    For the ``diagonal ideals'' of Example \ref{E:EgFrac1}, the intersection
	    of the diagonal with the Newton polygon is easily calculated 
	    using the dual vector \( v \).  The reader may check that its 
	    reciprocal is indeed 
	    \( v \cdot \one \).
	    (If it happens that \( v \cdot \one > 1\), then the log canonical
	    threshold is 1, and the multiplier ideal is trivial.)  See 
	    \cite{Kollar} for more details.
	\end{Ex}

	\begin{Ex} \label {E:complexlct}
	    To illustrate these ideas, we calculate the log-canonical 
	    threshold of a slightly more complicated ideal.  Let 
	    \[ {\mathfrak a} = (xy^{4}z^{6},x^{5}y,y^{7}z,x^{8}z^{8}).\]
	    After drawing the Newton polygon\footnote{Maple code illustrating
		this Newton polygon is available from the author by request.  Unfortunately,
		static \(2-\)dimensional representations are not very helpful.},
	    one sees that the diagonal
	    in \( {\mathbb R}^{3} \) intersects the triangular face generated by
	    the first three generators.  Therefore, the fourth generator 
	    \( x^{8}z^{8} \) can be ignored.  The intersection of the 
	    diagonal with the triangle whose vertices have coordinates
	    \( \{ (1,4,6), (5,1,0), (0,7,1) \} \) is the point 
	    \( (m,m,m) \), where \( m=\frac{191}{68} \).  The
	    log canonical threshold of \( {\mathfrak a} \) is \( \frac{1}{m} 
	    \), or \( \frac{68}{191} \).
	\end{Ex}
	
	The structure of the polygon \( P \) can in general be quite
	complicated, but it must have a single face which intersects
	the diagonal.  This face may not be simplicial, but it
	certainly decomposes into simplices, one of which intersects
	the diagonal in the same place and has no more than \( n \)
	vertices.  This demonstrates that the log canonical threshold
	of \( {\mathfrak a} \) is equal to that of a smaller ideal
	generated by no more monomials than the dimension \( n \) of
	the space.

	It has been conjectured \footnote{Actually, Shokurov's version of 
	this conjecture is stronger than that presented here.  It refers 
	to log canonical thresholds of effective Weil divisors on possibly 
	singular ambient spaces.} (\cite{Shokurov},\cite{Kollar}) that for 
	every dimension \( n \) the collection \( {\mathcal T}_{n} \) of all 
	log canonical thresholds satisfies the Ascending Chain Condition 
	(``All subsets have maximal elements'').  The restricted case of 
	ACC for monomial ideals follows from the fact that the partial 
	order of all monomial ideals has no infinite increasing sequences, nor 
	even any infinite antichains \cite{Maclagan}.  This fact doesn't 
	require any characterization of the thresholds.  If ACC is true, 
	then for any fixed dimension \( n \), there is a threshold \( 
	t_{n} \) closest to, but less than, one.  We attempted to use the 
	characterization above to calculate \( t_{n} \) in the monomial 
	case, but were unsuccessful.  It is known that \( t_{1} = 1/2, t_{2} = 
	5/6, t_{3}=41/42 \), Also, if we restrict to ideals of the form 
	\( {\mathfrak a} = (x_{1}^{b_{1}},\ldots,x_{n}^{b_{n}}) \), then it
	is known that we can 
	do no better than 
	\( t_{n} = \frac{a_{n}-1}{a_{n}} \), where \( a_{1}=2 \) and \( a_{n+1} = a_{n}^{2}+a_{n} \).
	The sequence \( a_{n} \) is \( (2, 6, 42, 1806, \ldots)\).
	We used a computer to calculate the log canonical threshold for large numbers of 
	monomial ideals, and found no evidence that the above pattern is 
	wrong in general.

\section*{Proof of the Theorem} \label {S:Mainproof}
	
	We will give a straightforward proof of the theorem, based on repeated
	blowups of the underlying space.  The basic proof structure is then 
	an induction, but this creates a problem:  After a single such 
	blowup \( f: {Y} \rightarrow {X} \) the space of interest is no longer 
	\( {{\mathbb A}^{n}}\), so an 
	inductive step doesn't apply.  
	
	This difficulty is not a serious one, 
	because \( Y \) is still locally \( {{\mathbb A}^{n}}\).  Also, all of the 
	above definitions can be extended to \( Y \) and onward.
	For example, the ``coordinate axes'' on \( Y \) should be taken to be
	proper transforms of those from \( X \), together with the exceptional 
	divisor(s).  The notion of a ``monomial ideal'' on \( X \) 
	generalizes on \( Y \) to an intersection of codimension-1 
	subschemes (monomials) supported on the ``coordinate axes.''
	These extensions are 
	consistent with those obtained by localizing on \( Y \) and
	identifying the coordinate patches with 
	\( {{\mathbb A}^{n}}\) in the obvious way.
	A briefer argument can be made 
	if one relies on the theory of toric varieties.  We will attempt to 
	point out these connections where appropriate.
	
	\begin{Def}
		By a {\bf monomial blowup}, we mean a blowup of \( X \) along 
		the intersection of some coordinate hyperplanes.  By a {\bf 
		sequence of monomial blowups}, we mean a sequence of blowups, 
		each of which is locally a monomial blowup.
	\end{Def}

	\begin{Def}
		Above we defined \( \one \) as a divisor on \( X = {{\mathbb A}^{n}} \),
		but we will need a more general notion.  If \( Y \) is obtained 
		from \( X \) by a sequence of monomial blowups, we let \( \one_{Y} \) be the 
		divisor which is the sum of the proper transforms of the coordinate 
		axes in \( X \), together which each exceptional divisor taken 
		with coefficient 1.  Thus \( \one \) is the union of the 
		``coordinate hyperplanes'' of \( Y \).  
		We regard \( \one_{Y} \) as an element of the 
		lattice \( {\mathbf L_{Y} }\), which must be defined as the free abelian 
		group on the coordinate hyperplanes in \( Y \).
	\end{Def}	
	
	The toric picture better illustrates what's 
	going on here:  the exceptional divisors and the 
	proper transforms of the coordinate axes 
	are precisely those effective divisors on \( Y \) which
	are invariant under the natural torus action.  Hence \( L_{Y} \)
	is the lattice of torically invariant divisors on \( Y \).
	In general, the sum of all of the effective toric divisors (each with
	coefficient one) on a toric variety is the
	anticanonical divisor.  So \( \one_{X} \) and 
	\( \one_{Y} \) are the torically natural anticanonical divisors, 
	and \( L_{Y} \) and \( L_{X} \) are the lattices of torically 
	invariant divisors.
	
	\begin{Lem} \label {L:oneminusone}
		Let \( X \) be \( {{\mathbb A}^{n}} \) or an 
		intermediate blowup, and let \( f: {Y} \rightarrow {X} \) 
		be a monomial blowup of 
		\( X \).  Then
		\[ \one_{Y} - f^{*}(\one_{X}) = K_{Y/X}.\]
	\end{Lem}	
	This can be seen without toric geometry by direct calculation; 
	it is easy to pull \( \one_{X} \) up to \( Y \) and count its multiplicity along 
	the exceptional divisor.
	
	\begin{Cor}
		If \( f: {Y} \rightarrow {X} \) arises by a sequence of 
		monomial blowups, then 
		\[ \one_{Y} - f^{*}(\one_{X}) = K_{Y/X}.\]
	\end{Cor}
	The corollary gives a convenient formula for \( K_{Y/X} \) when \( Y \) is a 
	log-resolution (via a sequence of monomial blowups) 
	of the monomial ideal \( {\mathfrak a} \).  It remains to see that such 
	a space \( Y \) exists:
	
	\begin{Lem} \label{L:ToricResolutionExists}
		Let \( X = {{\mathbb A}^{n}} \), and let \( {\mathfrak a} \) be a monomial 
		ideal on \( X \).  Then there is a sequence of monomial blowups 
		\( f: {Y} \rightarrow {X} \) which constitutes a log-resolution of 
		\( {\mathfrak a} \).
	\end{Lem}
	
	\begin{proof}
	    Here we must use some toric geometry.  The ideal \( {\mathfrak a} \)
	    defines a subset of the lattice \( L_{X} \).  The dual set of 
	    the ideal, \( \{ v \in L_{X}^{*} : \forall \lambda \in P,\ 
	    \langle v,\lambda \rangle \geq 1 \}\), defines a rational polytope \( P^{*} \) in the 
	    dual lattice \( L_{X}^{*} \).  To find a ``monomial log 
	    resolution'' of \( {\mathfrak a} \) is to find a sequence of toric 
	    blowups which refine the polytope \( P^{*} \) 
	    in the appropriate sense.  This
	    can be done because \( P^{*} \) is {\it rational}.  
	    The blowups required are exactly those required torically to 
	    resolve the singularity of the space \( Bl_{{\mathfrak a}}(X) \).
		Figure 2 indicates how this process might be used to resolve the cusp.
	    See \cite[section 2.6]{Fulton}, for more information on toric 
	    resolutions.
	\end{proof}

	\begin{figure}
		\hskip 1.4 in

		\mapleresult
		\begin{center}
		    \mapleplot{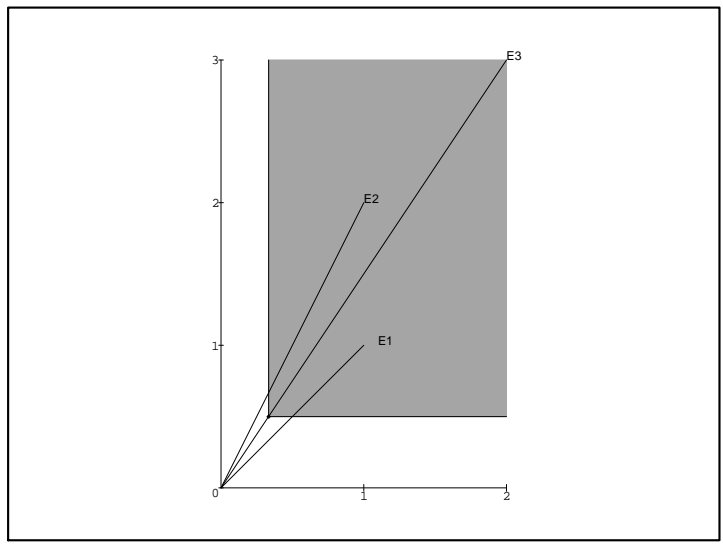}
		\end{center}
		\caption{The resolution of the cusp ideal \( (x^{3},y^{2}) \) 
		via dual polytope refinement.  The large
		rectangle is the dual polytope 
		\( P^{*} \) for this ideal.  The cusp is fully resolved because
		the ray representing \( E_{3}\) contains the lower left corner of the
		dual polytope, splitting it apart into sections without corners.}
	\end{figure}

	We now fix a monomial log-resolution \( f: {Y} \rightarrow {X} \),
	as in the Lemma.  We need to examine the relationship between 
	\( {\mathfrak a} \) and \( f^{-1}({\mathfrak a}) \).  By the definition of
	\( f \), \( f^{-1}({\mathfrak a}) \) is a line bundle.  It corresponds 
	to a divisor whose support is contained in the proper transforms of the coordinate 
	axes from \( X \) and the exceptional divisors.  We
	called the collection of such divisors \( L_{Y} \).  To 
	\( f^{-1}({\mathfrak a}) \) we may associate a single element 
	\( \gamma \) of 
	\( L_{Y} \), its ``generator.''  We may even give it a Newton polygon 
	\( P_{Y} \), namely the positive orthant translated to \( \gamma \).

	\begin{Lem} \label{L:Interior}
		Let \( f: {Y} \rightarrow {X} \) resolve \( {\mathfrak a} \)
		by a sequence of monomial blowups.  
		Let \( P_{X} \) be the Newton polygon of 
		\( {\mathfrak a} \), and let \( P_{Y} \) be as above.  
		Since \( f^{*} \) acts linearly on the lattices, we may 
		extend it to all of \( L_{X} \tensor {\mathbb R} \).  When we do this,
		
		1.  \( f^{*} \) takes the interior points of \( P_{X} \) to interior 
		points of \( P_{Y} \).
		
		2.  \( f^{*} \) takes the boundary points of \( P_{X} \) to boundary points
		of \( P_{Y} \).
		
		3.  \( f^{*} \) takes the points not in \( P_{X} \) to points not
		in \( P_{Y} \).
		
	\end{Lem}
	
	\begin{proof}
	    The lemma hinges on three basic ideas.  First, \( f^{*} \) is 
	    certainly a map from \( L_{X} \) to \( L_{Y} \), but because it
	    is linear, \( f^{*} \) extends to all of \( L_{X} \tensor 
	    {\mathbb R} \) in a natural way.  As a map of real vector spaces,
	    \( f^{*} \) is {\em continuous} because it is linear.	    
	    Second, for each of the effective toric divisors, or ``coordinate
	    planes'' \( E_{i} \) in \( Y \)
	    \[ord_{E_{i}}(f^{*}(\one_{X})) > 0.\]
	    The equality is strict because the blowups permitted are 
	    monomial.  
		Third, we have the standard equation
		\( f_{*}({\mathcal O}_{Y}(-E)) = \bar{\mathfrak a} \), 
		where \( \bar{\mathfrak a} \) is the integral 
		closure of \( {\mathfrak a} \).


		We will prove the lemma by proving part 3 first for integral points 
		\( \lambda \in L_{X} \tensor {\mathbb R} \), then for rational
		points, and finally for real points.  We will prove part 1 by using
		the strict
		positivity of \( f^{*}(\one_{X}) \).  Finally, we'll deduce part 2
		by continuity.
		
		Let \( \lambda \) be an integer point of \( L_{X} \) not in 
		\( P_{X} \).  Then 
		\( x^{\lambda} \notin \bar{\mathfrak a} = f_{*}({\mathcal O}_{Y}(-E)) \)
		so \( f^{*}(\lambda) \notin P_{Y} \).  If instead \( \lambda \notin 
		P_{X} \) has rational coordinates, then we can clear denominators.  
		Let \( n\lambda \) be integral.  
		\( n\lambda \notin nP({\mathfrak a}) = P({\mathfrak a}^{n}) \),
		so \( f^{*}(n\lambda) \notin f^{*}({\mathfrak a}^{n}) = 
		nP_{Y} \).  (Here we have used the just-proved integer case, as well
		as the fact that the resolution \( f: {Y} \rightarrow {X} \) resolving 
		\( {\mathfrak a} \) also resolves \( {\mathfrak a}^{n} \).)
		Dividing by \( n \) gives 
		\( f^{*}(\lambda) \notin P_{Y} \).  If 
		\( \lambda \notin P_{X} \) has real coordinates, choose a rational
		\( \mu \geq \lambda \) also not in \( P_{X} \).  
		\( f^{*}(\lambda) \leq f^{*}(\mu) \notin P_{Y} \), so 
		\( f^{*}(\lambda) \notin P_{Y} \).

		A standard convexity argument proves that if
		\( \lambda \in P_{X} \) then \( f^{*}(\lambda) \in P_{Y} \).
		To prove part 1 of the lemma, let \( \lambda \) be in the interior of
		\( P_{X} \).  Choose \( \mu \in P_{X} \) and \( \epsilon \in 
		{\mathbb R}^{+} \) with \(  \lambda = \mu + \epsilon \one  \).  Then
		\( f^{*}(\lambda) = f^{*}(\mu) + \epsilon f^{*}(\one) \).  
		\( f^{*}(\mu) \in P_{Y} \), and \( \epsilon f^{*}(\one) \) is strictly
		positive in every coordinate, so \( f^{*}(\lambda) \) is in the 
		interior of \( P_{Y} \).
		
		Part 2 of the lemma follows from the continuity of the map \( f^{*}\).
				
	\end{proof}

	
	We can now give the proof of the main theorem.  Because 
	\( \multr{{\mathfrak a}} \) is invariant under the natural torus action,
	it must be a monomial ideal.  We characterize the monomials 
	\( x^{\lambda} \) in \( \multr{{\mathfrak a}} \).  By definition, 
	\( x^{\lambda} \)
	is in \( \multr{{\mathfrak a}} \) if and only if 
	\[ div(f^{*}(x^{\lambda})) + K_{Y/X} - \floor{rE} \geq 0 \]
	(recall \( {\mathcal O}_{Y}(-E) = f^{-1}({\mathfrak a}) \)).  This condition simply 
	means that 
	\[ div(f^{*}(x^{\lambda})) + K_{Y/X} \text{ is in } \floor{rP_{Y}}\]
	(also recall \(\floor{rP_{Y}} =  \{ x^{\floor{\lambda}} : \lambda \in rP \}\)).
	Using the calculation of \( K_{Y/X} \) from Lemma \ref{L:oneminusone},
	this can be rewritten 
	\[ div(f^{*}(x^{\lambda})) - \one_{Y} + f^{*}(\one_{X}) 
		\in \floor{rP_{Y}}.\]
	This is of the form \( \{divisor\} - \one_{Y} \in \floor{rP_{Y}} \), 
	so we rewrite it as \( \{divisor\} \in int(rP_{Y}) \), obtaining
	\[ div(f^{*}(x^{\lambda})) + f^{*}(\one_{X}) 
		\in int(rP_{Y}).\]
	But this is just a condition on divisors from \( X \).  By Lemma 
	\ref{L:Interior}, parts 1 and 2, it is equivalent to
	\( (\lambda + \one_{X}) \in Int(rP_{X}) \).  The theorem is proved.


\bibliographystyle{plain}
\bibliography{monomials}
	
\end{document}